\documentclass[12pt,a4paper]{article}
\usepackage[utf8x]{inputenc}
\usepackage[english]{babel}
\usepackage[T1]{fontenc}
\usepackage{amsmath}
\usepackage{amsfonts}
\usepackage{amssymb}
\usepackage{makeidx}
\usepackage{graphicx}
\usepackage[left=3cm,right=2cm,top=3cm,bottom=2cm]{geometry}
\usepackage{setspace}
\usepackage{amsthm}
\usepackage{mathtools}
\usepackage{cite}
\usepackage[colorlinks=true,urlcolor=blue,
citecolor=red,linkcolor=blue,linktocpage,pdfpagelabels,
bookmarksnumbered,bookmarksopen]{hyperref}

\addto\captionsenglish{}

\newtheorem{definition}{Definition}[section]

\newtheorem{proposition}{Proposition}[section]

\newtheorem{theorem}{Theorem}[section]

\newtheorem{example}{Example}[section]

\newtheorem{lemma}{Lemma}[section]

\newtheorem{observation}{Remark}[section]

\newtheorem{corollary}{Corolary}[section]

\numberwithin{equation}{section}

\newcommand{\bo}{\begin{observation}}
\newcommand{\eo}{\end{observation}}
\newcommand{\bd}{\begin{definition}}
\newcommand{\ed}{\end{definition}}
\newcommand{\bp}{\begin{proposition}}
\newcommand{\ep}{\end{proposition}}
\newcommand{\bt}{\begin{theorem}}
\newcommand{\et}{\end{theorem}}
\newcommand{\bc}{\begin{corollary}}
\newcommand{\ec}{\end{corollary}}
\newcommand{\bl}{\begin{lemma}}
\newcommand{\el}{\end{lemma}}
\newcommand{\be}{\begin{example}}
\newcommand{\ee}{\end{example}}
\newcommand{\beq}{\begin{equation}}
\newcommand{\eeq}{\end{equation}}
\newcommand{\beqa}{\begin{equation*}}
\newcommand{\eeqa}{\end{equation*}}
\newcommand{\R}{\mathbb{R}}
\newcommand{\RN}{\mathbb{R}^{N}}
\newcommand{\N}{\mathbb{N}}
\newcommand{\Rdois}{\mathbb{R}^{2}}

\newcommand{\Ldois}{ L^{2}(\mathbb{R}^2) }
\newcommand{\Loito}{L^{\frac{8}{3}}(\Rdois)}

\newcommand{\B}{\mathcal{B}}
\newcommand{\T}{\mathcal{T}}
\newcommand{\A}{\mathcal{A}}

\newcommand{\Hh}{\mathcal{H}}

\newcommand{\Ss}{\mathcal{S}}
\newcommand{\m}{\mathfrak{m}}
\newcommand{\Aa}{\mathfrak{a}}
\newcommand{\Bb}{\mathfrak{b}}
\newcommand{\Dd}{\mathfrak{d}}

\newcommand{\M}{\mathcal{M}}
\newcommand{\K}{\mathcal{K}}
\newcommand{\I}{\mathcal{I}}

\newcommand{\Hum}{H^{1}(\mathbb{R}^2)}
\newcommand{\Ls}{ L^{s}(\mathbb{R}^2)}

\newcommand{\intR}{\displaystyle\int\limits_{\mathbb{R}^2}}
\newcommand{\un}{u_{n}}

\newcommand{\vn}{v_{n}}

\newcommand{\yn}{y_{n}}
\newcommand{\wn}{w_{n}}

\newcommand{\until}{\tilde{u}_n}
\newcommand{\CF}{\rightharpoonup}

\newcommand{\RA}{\rightarrow}

\newcommand{\ds}{\displaystyle\int\limits}
\newcommand{\intlog}{ \displaystyle\int\limits_{\mathbb{R}^2} \displaystyle\int\limits_{\mathbb{R}^2} \ln (|x-y|)u^2(x)u^2(y) dx dy}

\usepackage{times, color, xcolor}
\addto\captionsportuguese{
}

\begin{document}

 	\title{Existence of normalized solutions for the planar Schr\"odinger-Poisson system with exponential critical nonlinearity
		}
	\author{
		Claudianor O. Alves$^{1}$ \thanks{Corresponding author} \footnote{E-mail address: coalves@mat.ufcg.edu.br }, Eduardo  de S. Böer$^{2}$  \thanks{ E-mail address: eduardoboer04@gmail.com, Tel. +55.51.993673377} \,\, and \,\, Ol\'{\i}mpio H. Miyagaki$^{2}$ \footnote{ E-mail address: ohmiyagaki@gmail.com, Tel.: +55.16.33519178 (UFSCar).}  \\
		{\footnotesize $^{1}$  Unidade Acad\^emica de Matem\'atica, Universidade Federal de Campina Grande}\\ {\footnotesize 58109-970 Campina Grande-PB-Brazil}\\
		{\footnotesize $^{2}$  Department of Mathematics, Federal University of S\~ao Carlos,}\\
		{\footnotesize 13565-905 S\~ao Carlos, SP - Brazil}  
		}
\noindent
				
	\maketitle

\noindent \textbf{Abstract:} In the present work we are concerned with the existence of normalized solutions to the  following Schrödinger-Poisson System
$$
\left\{ \begin{array}{ll}
-\Delta u + \lambda u + \mu (\ln|\cdot|\ast |u|^{2})u = f(u) \textrm{ \ in \ } \mathbb{R}^2 , \\
\intR |u(x)|^2 dx = c,\ c> 0  ,
\end{array}
\right. 
$$ 
for $\mu \in \R $ and a nonlinearity $f$ with exponential critical growth. Here $ \lambda\in \R$ stands as a Lagrange multiplier and it is part of the unknown. Our main results extend and/or complement some results found in \cite{Ji} and \cite{[cjj]}.

\vspace{0.5 cm}


\section{Introduction}

In the present paper we are interested with the existence of normalized solutions for the following Schrödinger-Poisson System,
\beq\label{i5}
\left\{ \begin{array}{rclcl}
i\psi_t -\Delta\psi + \tilde{W}(x)\psi + \gamma \omega \psi &= &0 \  &  \textrm{ in \ } &  \RN\times \R \\
\Delta \omega &=& |\psi|^2 & \textrm{ in \ } & \RN,
\end{array} \right. 
\eeq
where $ \psi : \RN \times \R \RA \mathbb{C} $ is the time-dependent wave function, $ \tilde{W}: \RN \RA \R $ is a real external potential and $ \gamma > 0 $ is a parameter. The function $ \omega $ stands as an internal potential for a nonlocal self-interaction of the wave function $ \psi $. The usual ansatz $ \psi(x, t)=e^{-i\Theta t}u(x) $, with $ \Theta \in \R $, for standing wave solutions of (\ref{i5}) leads to
\beq\label{i1}
\left\{ \begin{array}{rclcl}
-\Delta u + W(x)u + \gamma \omega u & = &0  & \textrm{ in \ }& \RN \\
\Delta \omega &= &u^2 &\textrm{ in \ }& \RN,
\end{array} \right. 
\eeq
with $ W(x)= \tilde{W}(x) + \Theta $. From the second equation of (\ref{i1}) we observe that $ \omega : \RN \RA \R $ is determined only up to harmonic functions. In this point of view, it is natural to choose $ \omega $ as the Newton potential of $ u^2 $, i.e., $ \Gamma_N \ast u^2 $, where $ \Gamma_N $ is the well-known fundamental solution of the Laplacian
$$
\Gamma_N (x) = \left\{ \begin{array}{rcl}
\dfrac{1}{N(2-N)\sigma_N}|x|^{2-N} &\textrm{ if \ } &N\geq 3 ,\\
\dfrac{1}{2\pi}\ln |x|& \textrm{ if \ }& N= 2 ,
\end{array} \right.
$$
in which $ \sigma_N $ denotes the volume of the unit ball in $ \RN $. From this formal inversion of the second equation in (\ref{i1}), as it is detailed in \cite{[4]}, we obtain the following integro-differential equation
\beq\label{i2} 
-\Delta u + W(x) u + \gamma ( \Gamma_N \ast |u|^2) u = b|u|^{p-2}u,\ p> 2, \ b>0, \textrm{ \  \ in \ } \RN.
\eeq

Then, we make a quick overview of the literature. To begin with, we note that the case $ N=3 $ has been extensively studied, due to its relevance in physics. It is curious that, although this equation is called ``Choquard equation'', it has first studied by Fröhlich and Pekar in \cite{[12] , [11] , [22]}, to describe the quantum mechanics of a polaron at rest in the particular case $ W(x)\equiv a > 0 $ and $ \gamma > 0 $. Then, it was introduced by Choquard in 1976, to study  an electron trapped in its hole. From the applied point of view, the local nonlinear term on the right side of equation (\ref{i2}) usually appears in Schrödinger equations as a model to the interaction among particles. 
%For example, Penrose has derived equation (\ref{i2}) while discussing about the self gravitational collapse of a quantum-mechanical system in  \cite{[18]}. Besides that, many results have been derived recently about the existence and regularity for this kind of problem, including for more general convolution potentials, such as in  the evolution equation $ i\partial_t \phi = \Delta \phi +(V \ast |\phi|^2)\phi,  $ which models the interaction of a large system of nonrelativistic bosonic atoms and molecules.  

We could discuss a bunch of variations of these kind of equations in the three dimensional case, but in order to make it concise, we refer the readers to the following papers, and the references therein, \cite{[14] , [16] , [1] , [2] , [3] , [20] , [18]}. 

Once we turn our attention to the case $ N=2 $, we immediately see that the literature is scarce. In this case, we can cite the recent works of \cite{[6], [cjj], [10], [alves], [boer], [boer3]}. In \cite{[6], [10]}, the authors have proved the existence of infinitely many geometrically distinct solutions and a ground state solution, considering $ W $ continuous and $ \mathbb{Z}^2 $-periodic and  $ W(x)\equiv a > 0 $, respectively, and the particular case $ f(u)=b|u|^{p-2}u $. Then, in \cite{[alves]}, the authors have dealt with equation (\ref{P}) considering a nonlinearity with exponential critical growth. They have proved the existence of a ground state solution via minimization over Nehari manifold. Moreover, in \cite{[boer3]}, the authors have proved existence and multiplicity results for the $ p-$fractional Laplacian operator. Finally, in \cite{[boer]}, the authors deal with a $ (p, N)$-Choquard equation and prove existence and multiplicity results.

We call attention to \cite{[cjj]} in which work the authors have dealt with the existence of stationary waves with prescribed norm considering $ \lambda \in \R $ and that consists in a key reference to our work. Another important reference is \cite{[olimpio]}, since we manage to adapt some techniques from both of them. We also refer to \cite{[olimpio]} for a relevant overview about Schr\"odinger prescribe norm problems.

 % Also, it is interesting to mention that the application of variational methods in the planar case 
%\beq\label{i3}
%-\Delta u + a(x) u + \dfrac{\gamma}{2\pi} ( \ln(|\cdot|) \ast |u|^2) u = 0 \textrm{ \  \ in \ } \Rdois ,
%\eeq
%took place after the work of Stubbe \cite{[21]}, in which he set up a variational framework for (\ref{i3}) within a subspace of $ \Hum $, where the associated functional is well-defined. 

In the present paper we focus on finding prescribe norm solutions for the planar equation 
\begin{equation} \label{P}
\left\{ \begin{array}{ll}
-\Delta u + \lambda u + \mu (\ln|\cdot|\ast |u|^{2})u = f(u) \textrm{ \ in \ } \mathbb{R}^2 \\
\intR |u(x)|^2 dx = c, c> 0 
\end{array}
\right. .
\end{equation}
where $ \lambda, \mu \in \R $ and $f: \R \RA \R $ is continuous, with primitive $ F(t)=\int\limits_{0}^{t}f(s)ds $.  This approach seems to be particularly meaningful from the physical point of view, because there is a conservation of mass.

The main difficulties in the proof our main results are associated with the fact that we are working with critical nonlinearities in the whole $\mathbb{R}^2$ and with the logarithm term, which are unbounded and changes sign.

Our aim is to extend and/or complement the results already obtained in the literature and cited above, more precisely the results found in \cite{Ji} and \cite{[cjj]}, by working with a nonlinearity that has an exponential critical growth. We recall that a function $ h $ has an exponential subcritical growth at $ +\infty $, if
$$
\lim\limits_{s\RA + \infty}\dfrac{h(t)}{e^{\alpha t^2}-1} = 0 \textrm{ \ , for all \ } \alpha >0 ,
$$
and we say that $ h $ has $ \alpha_0 $-\textit{critical} exponential growth at $ +\infty $, if
$$
\lim\limits_{t\RA + \infty}\dfrac{h(t)}{e^{\alpha t^2}-1} = \left\{ \begin{array}{ll}
0, \ \ \ \forall \ \alpha > \alpha_{0}, \\
+\infty , \ \ \ \forall \ \alpha < \alpha_0.
\end{array} \right. 
$$
As usual conditions while dealing with this kind of growth, found in works such as \cite{[19] , [9]}, we assume that $ f $ satisfies
$$f\in C(\R , \R), f(0)=0 \mbox{ and has a critical exponential growth with } \alpha_0 = 4\pi . \leqno{(f_1)}$$
$$ \lim\limits_{|t|\RA 0} \dfrac{|f(t)|}{|t|^\tau }=0 \mbox{, for some } \tau > 3.  \leqno{(f_2)}$$
From $ (f_1) $ and $(f_2)$, given $ \varepsilon >0 $, $\alpha > \alpha_0,$ fixed, for all  $ p>2 $, we can find two constants $ b_1=b_1(p, \alpha , \varepsilon) > 0 $ and $ b_2=b_2(p, \alpha , \varepsilon) > 0$ such that 
\begin{equation}\label{eq1}
f(t)\leq \varepsilon |t|^\tau +b_1 |t|^{p-1}(e^{\alpha t^2}-1) \ , \ \ \ \forall \ t \in \R,
\end{equation} 
and
\begin{equation}\label{eq2}
F(t) \leq \varepsilon|t|^{\tau + 1} + b_2 |t|^p(e^{\alpha t^2}-1) \ , \ \ \ \forall t \in \R .
\end{equation}
In order to verify that $(PS)$ sequences are bounded in $\Hum$, we will need the following conditions: 
$$ \mbox{ there exists }\ \theta > 6 \ \mbox{ such that}\   f(t)t \geq \theta F(t) > 0, \ \mbox{for all }\  t\in \R \setminus \{0\}, \leqno{(f_3)}$$
$$ \mbox{ there  exist}\  q>4 \ \mbox{ and} \  \nu > \nu_0 \ \mbox{ such that}\  F(t) \geq \nu |t|^q , \ \mbox{for all} \  t \in \R . \leqno{(f_4)}$$

Next we provide some important definitions for our work. Once we will use variational techniques, we consider the associated Euler-Lagrange functional $ I: \Hum \RA \R \cup \{ \infty\} $ given by
\begin{equation}\label{I}
I(u) = \dfrac{1}{2}A(u) + \dfrac{\mu}{4}V(u) -\intR F(u) dx, 
\end{equation}
where
$$
A(u) = \intR |\nabla u|^2 dx \mbox{ \ \ \ and \ \ \ } V(u)= \intlog  .
$$
It is easy to verify, from Moser-Trudinger inequality, Lemma \ref{l5}, and Hardy-Littlewood-Sobolev inequality (HLS) (found in \cite{[15]}), that $ I $ is well-defined on the slightly smaller Hilbert space
\beq \label{X}
X = \left\{ u\in \Hum ; ||u||_{\ast}^{2}=\intR \ln(1+|x|)u^2(x) dx < \infty \right\} \subset \Hum ,
\eeq
endowed with the norm $ ||\cdot||_{X}=\sqrt{||\cdot||^2 + ||\cdot||_{\ast}^{2}} $, where $ ||\cdot|| $ is the usual norm in $ \Hum $. Moreover, $I$ is $ C^1 $ on $ X $ (see \cite{[alves]} and \cite{[6]}) and any critical point $ u $ of $ I\big\vert_{S(c)} $ corresponds to a solution of (\ref{P}), where $ \lambda \in \R $ appears as a Lagrange multiplier and
$$
S(c) = \{ u \in X \ ; \ ||u||_{2}^{2}=c \} .
$$

The two first results of the present paper involving the existence of solution are the following:

\bt \label{t1}
Suppose that $ f $ satisfies $ (f_1)-(f_4) $. Then, there are $\mu_0,\nu_0>0$ such that  problem (\ref{P}) has at least one weak solution $ u\in S(c) $ such that $ I(u) = \m_{\nu} > 0 $ for $ \mu \in (0, \mu_0) $, $ c\in (0, 1) $ and $ \nu > \nu_0 $, where
\begin{equation}\label{eq10}
\m_\nu = \inf\limits_{\gamma \in \Gamma} \max\limits_{t \in [0,1]}I(\gamma(t)), 
\end{equation} 
with $ \Gamma = \{ \gamma \in C([0, 1], S(c)) \ ; \ \gamma(0)=u_1 \mbox{ and } \gamma(1)=u_2\} $, $ u_1, u_2 \in S(c)$. 
\et

In a similar way, we get the following result.

\bt \label{t2}
Suppose that $ f $ satisfies $ (f_1)-(f_4) $ and $ \mu > 0 $. Then, there are $c_0,\nu_0$ such that problem (\ref{P}) has at least one weak solution $ u\in S(c) $ such that $ I(u) = \m_{\nu} > 0 $ for  $ c \in (0, c_0) $ and $ \nu > \nu_0 $.  
\et

Once the proof of Theorem \ref{t2} is similar to the proof of Theorem \ref{t1}, where the only change is  the control of the inequalities from parameter $ \mu $ to the mass $ c $, we will not write it down.

As an immediate consequence from Theorems \ref{t1} and \ref{t2} we get the following type of least energy level. We do not call it a ground state since it cannot be take over all the possible solutions for \eqref{P}.

\bc
Under the hypotheses of Theorem \ref{t1} (or Theorem \ref{t2}) problem \eqref{P} has a least energy level solution $ u \in S(c) $, in the sense that, there is a function $ u \in S(c) $ satisfying
$$
I(u) = \m_l = \inf \{ I(v) \ ; \ I\big\vert_{S(c)}'(v) =0 \mbox{ \ and \ } Q(v) = 0 \}, 
$$
where $ Q $ is defined in \eqref{Q}.
\ec

Related to the existence of multiple solutions, we will use a genus approach that is based on the ideas of \cite{Ji}. In order to do so, we need the following condition.
$$f\in C(\R , \R), f(0)=0, f \mbox{ is odd and has a critical exponential growth with } \alpha_0 = 4\pi . \leqno{(f_1 ')}$$

As in the existence case, we have two results which differ by the way that we control some estimates, either by the parameter $ \mu $ or by the mass $ c $.

\bt\label{t3}
Suppose that $ f $ satisfies $ (f_1 ')-(f_4) $, $ c\in (0, 1) $ and $ \mu\in (0, \mu_1) $, for $\mu_1$ defined in \eqref{eqe1}. Then, given $ n \in \N $, there is $ \tilde{\nu} =\tilde{\nu}(n)>0$ sufficiently large such that, if $ \nu \geq \tilde{\nu} $, \eqref{P} has at least $ n $ non-trivial weak solutions , $ u_j \in S(c) $, verifying $ I(u_j) < 0 $, for $ 1 \leq j \leq n $.
\et

\bt\label{t4}
Suppose that $ f $ satisfies $ (f_1 ')-(f_4) $, $ c\in (0, \min\{c_1, 1\}) $, , for $c_1$ defined in \eqref{eqe1}, and $ \mu \in \R^+ $. Then, given $ n \in \N $, there is $ \tilde{\nu} =\tilde{\nu}(n)>0 $ sufficiently large such that, if $ \nu \geq \tilde{\nu} $, \eqref{P} has at least $ n $ non-trivial weak solutions , $ u_j \in S(c) $, verifying $ I(u_j) < 0 $, for $ 1 \leq j \leq n $.
\et

Before concluding this section, we would like point out that our main results complement the study made in \cite{[cjj]} because in that paper it was not considered the case where the nonlinearity has a critical exponential growth, while in \cite{Ji}, the authors neither consider the existence of normalized solutions nor the presence of an unbounded indefinite internal potential.

The paper is organized as follows: in Section 2 we present some technical and essential results, some of them already derived in previous works. Section 3 is devoted to the study of the geometry of the associated functional and some convergence results. Section 4 consists in the proof of our existence main results and, finally, Section 5 is concerned with the proof of the multiplicity results.

Throughout the paper, we will use the following notations:
\begin{itemize}
	
	\item  We fix the values $ r_1, r_2 >0 $ such that $ r_1> 1 $, $ r_1 \sim 1 $ and $ \frac{1}{r_1}+\frac{1}{r_2}=1 $.
	
	\item $ \Ls $ denotes the usual Lebesgue space with norm $ ||\cdot ||_s $.
	
	\item $ X' $ denotes the dual space of $ X $. 
	
	\item   $ B_r(x) $ is the ball centered in $ x $ with radius $ r>0 $, simply $ B_r $ if $ x=0 $.
	
	\item   $ K_i, b_i $, $ i\in \N $, stand for important constants that appear in the estimates obtained.
	
	\item  $ C_i $, $ i\in \N $, will denote different positive constants whose exact values are not essential to the exposition of arguments.

\end{itemize}

\section{Framework and some Technical Results}

In this section, we will focus in presenting additional framework properties and a few technical results. Some of them can be found in \cite{[6], [cjj], [olimpio], [boer]} and we will omit their proofs here.  
 
We begin defining three auxiliary symmetric bilinear forms 
\beqa
(u , v) \mapsto B_1(u, v)=\intR \intR \ln(1+|x-y|)u(x)v(y) dx dy ,
\eeqa
\beqa
(u , v) \mapsto B_2(u, v)=\intR \intR \ln\left(1+\dfrac{1}{|x-y|}\right)u(x)v(y) dx dy ,
\eeqa
\beqa
(u , v) \mapsto B(u, v)=B_1(u, v)-B_2(u, v)=\intR \intR \ln(|x-y|)u(x)v(y) dx dy. 
\eeqa
The above definitions are understood to being over measurable functions $u, v: \Rdois \RA \R $, such that the integrals are defined in the Lebesgue sense. Then, we can define the functionals  $V_1: \Hum \RA [0, \infty],$ $V_2: L^{\frac{8}{3}}(\Rdois) \RA [0, \infty) $ given by $V_1(u)=B_1(u^2 , u^2)$ and $V_2(u)=B_2(u^2 , u^2) $, respectively. Moreover, one should observe that $V(u)=V_1(u)-V_2(u)$.

\bo\label{obs1}
\begin{itemize}
\item[(i)] From $ 0 \leq \ln(1+r)\leq r $, for $ r>0 $, and (HLS), for $ u, v \in L^{\frac{4}{3}}(\Rdois) $, 
$|B_2(u, v)|\leq K_0 ||u||_{\frac{4}{3}}||v||_{\frac{4}{3}}$ , where   $ K_{0}>0 $ is the (HLS) best constant. Hence, 
\beq \label{v2}
|V_2(u)|\leq K_0||u||_{\frac{8}{3}}^{4} \ , \ \ \forall \ u\in L^{\frac{8}{3}}(\Rdois) ,
\eeq

\item[(ii)] We recall the so-called Gagliardo-Nirenberg inequality
\begin{equation}\label{eq3}
||u||_s \leq K_{GN}^{\frac{1}{s}}||\nabla u||_{2}^{\sigma}||u||_{2}^{1-\sigma} , \mbox{ for \ } t> 2 \mbox{ \ and \ } \sigma = 2\left(\dfrac{1}{2}-\dfrac{1}{s}\right) .
\end{equation} 

\item[(iii)] From equations (\ref{v2}) and (\ref{eq3}), we obtain a positive constant $ K_1 > 0 $ such that
\begin{equation}\label{eq4}
|V_2(u)| \leq K_1 c^{\frac{3}{2}} \sqrt{A(u)} , \forall \ u \in \Hum . 
\end{equation}

\item[(iv)] Observing that
\beq \label{eq5}
\ln(1+|x-y|) \leq \ln(1+|x|+|y|) \leq \ln (1+|x|)+ \ln(1+|y|), \textrm{ \ for all \ } x, y\in \Rdois,
\eeq
we obtain the estimate
\beq\label{eq6}
B_1(uv, wz) \leq ||u||_{\ast}||v||_{\ast}||w||_2 ||z||_2 + ||u||_2 ||v||_2 ||w||_{\ast}||z||_{\ast} ,
\eeq
for all $ u, v, w, z \in \Ldois $. \\

\end{itemize}
\eo

We are going to need the following results from \cite{[6]}.

\bl \label{l1} (\cite[Lemma 2.2]{[6]})
(i) The space $ X $ is compactly embedded in $ \Ls $, for all $ s\in [2, \infty) $. \\
(ii) The functionals $ V_0, V_1, V_2 $ and $ I $ are of class $ C^1 $ on $ X $. Moreover, $ V_{i}'(u)(v)=4B_i(u^2, uv) $, for $ u, v \in X $ and $ i=0, 1, 2 $. \\
(iii) $ V_2 $ is continuous (in fact continuously differentiable) on $\Loito$ .
\el

\bl \label{l2}
(\cite[Lemma 2.1]{[6]}) Let $(\un)$ be a sequence in $\Ldois $ and $ u\in \Ldois \setminus \{0\} $ such that $\un \RA u $ pointwise a.e. on $\Rdois $. Moreover, let $(\vn)$ be a bounded sequence in $\Ldois$ such that
$$
\sup\limits_{n\in \N}B_1(\un^2, \vn^2) < \infty .
$$
Then, there exist $n_0\in \N$ and $ C>0 $ such that $ ||\un||_{\ast}< C $, for $ n\geq n_0 $. If, moreover, 
$$
B_1(\un^2, \vn^2)\RA 0 \textrm{ \ \ and \ \ } ||\vn||_{2}\RA 0, \textrm{ \ as \ } n\RA \infty ,
$$
then 
$$
||\vn||_{\ast}\RA 0 \textrm{ \ ,  as \ } n \RA \infty .
$$
\el

\bl\label{l3}
(\cite[Lemma 2.6]{[6]}) Let $ (\un) $, $ (\vn) $ and $ (\wn) $ be bounded sequences in $ X $ such that $ \un \rightharpoonup u $ in $ X $. Then, for every $ z\in X $, we have $ B_1 (\vn \wn \ , \ z(\un - u))\RA 0 $, as $ n\RA + \infty $.
\el

Also, we borrow the following result from \cite{[cjj]}.

\bl\label{l4}
(\cite[Lemma 2.5]{[cjj]}) Let $ (\un)\subset S(c) $ and assume the existence of $ \varepsilon \in (0, c) $ such that for all $ R>0 $, we have 
$$
\liminf_{n} \sup\limits_{y\in \mathbb{Z}^2} \ds_{B_R(y)} |\un|^2\,dx \leq c- \varepsilon .
$$
Then, 
$$
\limsup_{n} V_1(\un) = + \infty .
$$
\el

\bc\label{c2}
Let $ (\un) \subset S(c) $ and assume that $ (V_1(\un)) $ is bounded. Then, there exists $ R_0 \geq 2 $ such that 
$$
\liminf_{n} \sup\limits_{y\in \mathbb{Z}^2} \ds_{B_{R_{0}}(y)} |\un|^2\, dx > \frac{c}{2} .
$$
\ec
\begin{proof}
From Lemma \ref{l4}, for $ \varepsilon = \frac{c}{2} $, there exists $ \tilde{R} > 0 $ satisfying 
$$
\liminf_{n} \sup\limits_{y\in \mathbb{Z}^2} \ds_{B_{\tilde{R}}(y)} |\un|^2 \,dx > \frac{c}{2} .
$$
Consequently, the corollary follows considering $ R_0 = \tilde{R} + 2 $.
\end{proof}

Now, we turn our attention to the term with exponential critical growth. First of all, we remember the well-known Moser-Trudinger inequality.

\bl\label{l5}
\cite{[5]} If $\alpha >0$ and $ u\in \Hum $, then
$$
\intR \left(e^{\alpha |u|^2} - 1 \right)\,dx  < +\infty .
$$
Moreover, if $ ||\nabla u||_{2}^{2}\leq 1 $, $ ||u||_{2}^{2}\leq M < \infty $ and $ \alpha <  4 \pi $, then there exists $ K_{\alpha, M}=C(M, \alpha) $, such that
$$
\intR \left(e^{\alpha |u|^2} - 1 \right) dx < K_{\alpha, M} .
$$
\el

Then, inspired by \cite{[olimpio]}, we prove the following corollary.

\bc\label{c1}
Let $ (\un)\subset S(c) $ satisfying
$
\limsup A(\un) <  1 - c .
$ Then, for all $ p \geq 1 $, there exist values $ \beta > 1 $, $ \beta \sim 1 $, $ p_0 = p_0(p) > 2 $ and a constant $ K_2 = K_2 ( \beta, c) >0 $ such that
$$
\intR |\un|^p (e^{4\beta \pi |\un(x)|^2} - 1) dx \leq K_2 ||\un||_{p_0}^{p} , \forall \ n \in \N .
$$ 
\ec
\begin{proof}
Since $\displaystyle \limsup_{n} A(\un) < 1 -c $, there exist $ d>0 $ and $ n_0 \in \N $ such that $ ||\un||^2 < d < 1 $, for all $ n \geq n_0 $. Thus, there exists $ \beta >1 $, $ \beta \sim 1 $ with $ 1< \beta < \frac{1}{d} $ and $ s_1 > 1 $, $ s_1 \sim 1 $ with $ 1 < s_1 < \frac{1}{d \beta} $. Let $ s_2 > 2 $ such that $ \frac{1}{s_1}+\frac{1}{s_2}=1 $. Consequently, 
$$
\intR |\un|^p (e^{4\beta \pi |\un(x)|^2} - 1) dx \leq ||\un||_{ps_2}^{p}\left(\intR (e^{4\beta s_1 d \pi \left(\frac{\un(x)}{||\un||}\right)^2} - 1) dx\right)^{\frac{1}{s_1}} \leq C_1 ||\un||_{ps_2}^{p} ,
$$
for all $ n \geq n_0 $. Now, consider
$$
C_2 = \max\left\{ \left(\intR (e^{4\beta s_1 \pi |u_1 (x)|^2} - 1) dx\right)^{\frac{1}{s_1}} , ... \ , \left(\intR (e^{4\beta s_1 \pi |u_{n_0} (x)|^2} - 1) dx\right)^{\frac{1}{s_1}}\right\} .
$$
Then, for each $ 1\leq n \leq n_0 $, we have
$$
\intR |\un|^p (e^{4\beta \pi |\un(x)|^2} - 1) dx \leq C_2 ||\un||_{ps_2}^{p} .
$$
Therefore, the lemma follows setting $ p_0 = p s_2 > 2 $ and $ K_2 = \max\{ C_1 , C_2\}>0 $.
\end{proof}

From now on, unless we say otherwise, $ \beta >0 $ will stand as that given by Corollary \ref{c1}. The next lemma plays an important role in our work but since its proof is very similar to that in \cite[Lemma 3.7]{[boer3]}, we omit it here.

\bl\label{l6}
Let $ (\un)\subset X $ be bounded in $ \Hum $ and $ R > 0  $ satisfying
\begin{equation}\label{eq7}
\liminf_{n} \sup\limits_{y\in \mathbb{Z}^2} \ds_{B_R(y)}\un^2(x) dx > 0 .
\end{equation}
Then, there exists $ u\in \Hum \setminus\{0\} $ and a sequence $ (\yn)\subset \mathbb{Z}^2 $ such that, up to a subsequence, $ \tilde{u}_n(x)=u_n(\cdot-y_n) \rightharpoonup u$ in $ \Hum $. 
\el

In the last remark of this section we discuss the geometry of a real function that appears in some further arguments.

\bo\label{obs5}
Consider the real function $ h: (0, + \infty) \RA (0, + \infty) $ given by $ h(t)=\Aa t-\Bb t^{\frac{1}{2}}+\Dd $, for $ \Aa , \Bb , \Dd > 0 $. Note that $ h'(t)=\Aa - \frac{\Bb}{2t^{\frac{1}{2}}} $ and $ h''(t)=\frac{\Bb}{4t^{\frac{3}{2}}} $. One can see that $ h'(t) = 0 $ for $ t = \left(\frac{\Bb}{2\Aa}\right)^{2} $ and $ h''(t) > 0 $ for all $ t \in (0, + \infty) $. Consequently, $ h $ is convex and $ \overline{t} = \left(\frac{\Bb}{2\Aa}\right)^{2} $ is a global minimum for $ h $. Moreover, $ h(\overline{t}) = -\frac{\Bb^2}{4\Aa}+ \Dd $. Thus, we conclude that
$$
h(\overline{t}) \geq 0 \mbox{ \ \ if and only if \ \ } \dfrac{\Bb^2}{4\Aa} \leq \Dd .
$$
\eo 

\section{Geometry of $ I $ and key auxiliary results}

The present section will be devoted to derive some geometrical properties of the functional $ I $ and some very important auxiliary results, such as those involving boundedness and convergence of sequences on $ X $. We start defining the map $ H : \Hum \RA \R $ by $ H(u, t) = e^t u(e^t x)  $, for all $ x\in \Rdois $ and $ t\in \R\setminus \{0\} $, and, for each $ u\in \Hum $ fixed, the function $ \varphi_u(t): \R \RA \R $ given by $ \varphi_u(t) = I(H(u, t)) $.

\bo\label{obs3}
For any $ u\in S(c) $, one should observe that $ H(u, t) \in S(c) $, for all $ t\in \R\setminus\{0\} $. Moreover, we have the following:
\begin{itemize}
\item[(i)] $ V(H(u, t))=V(u) - t ||u||_{2}^{4} $, 

\item[(ii)] $ A(H(u, t ))=e^{2t}A(u) $, 

\item[(iii)] $ ||H(u, t)||_{p}^{p} = e^{(p-2)t}||u||_{p}^{p} $, for all $ p \geq 1 $.
\end{itemize}
\eo

\bl\label{l7}
Assume $ (f_1)-(f_4) $ and $ u\in S(c) $. Then, 
\begin{itemize}
\item[(1)] $ A(H(u, t))\RA + \infty $ and $ \varphi_u(t) \RA - \infty $, as $ t\RA + \infty $ .

\item[(2)] $ A(H(u, t))\RA 0 $ and 
$
\left\{ \begin{array}{ll}
\varphi_u(t) \RA + \infty \mbox{ , if \ } \mu > 0 \\
\varphi_u(t) \RA - \infty \mbox{ , if \ } \mu < 0 
\end{array} \right. \mbox{ , as \ } t \RA - \infty .
$

\item[(3)] $ \varphi_u(t) \RA I(u) $, as $ t \RA 0 $. 
\end{itemize}
\el
\begin{proof}
\textbf{(1)} From Remark \ref{obs3}-(ii), $ A(H(u, t))\RA + \infty $, as $ t\RA + \infty $. Moreover, from condition $ (f_4) $, Remark \ref{obs3} and $ q>4 $, we have
$$
\varphi_u(t) \leq \dfrac{e^{2t}}{2}A(u) + \dfrac{\mu}{4}V(u) - \dfrac{\mu}{4}c^2 t - \nu e^{(q-2)t}||u||_{q}^{q}\RA - \infty \mbox{ , as } t \RA + \infty .
$$
\textbf{(2)} Once again, from Remark \ref{obs3}, $ A(H(u, t))\RA 0 $ and $ ||H(u, t)||_{p}^{p}\RA 0 $, as $ t\RA - \infty $, for all $ p > 2 $. Thus, there are $ t_0 < 0 $ and $ d\in (0, 1) $ such that $ \beta d r_1 < 1 $ and $ ||H(u, t)||^2 \geq d $, for all $ t \in (- \infty , t_0] $, with $ r_1 > 1 $, $ r_1 \sim 1 $ and $ \frac{1}{r_1}+\frac{1}{r_2}=1 $, $ \beta > 1 $, $ \beta \sim 1 $.

Then, from (\ref{eq2}), we have
$$
|F(H(u, t))| \leq \varepsilon |H(u, t)|^{\tau +1 }+b_2 |H(u, t)|^{q}(e^{4\pi \beta |H(u, t)|^2}-1) , \ \forall \ t \leq t_0 ,
$$ 
and applying Hölder inequality and Lemma \ref{l5}, we obtain a constant $ C_1 >0 $ such that
$$
\intR F(H(u, t)) dx \leq \varepsilon||H(u, t)||_{\tau + 1}^{\tau + 1} + C_1 ||H(u, t)||_{qr_2}^{q}\RA 0 \mbox{ , as \ } t\RA - \infty .
$$
Consequently, from this fact and from Remark \ref{obs3}, we get item (2).
\end{proof}   

For the next result, let us consider the subsets of $ S(c) $: 
$$
\A^{-}= \{ u \in S(c) \ ; \ V(u) < 0 \} ,  \A^{+}= \{ u \in S(c) \ ; \ V(u) \geq 0 \}, \A_{r}= \{ u \in S(c) \ ; \ A(u) \leq r \} ,
$$
$ \A_{r}^{+} = \A^{+} \cap \A_r $ and $ \A_{r}^{-} = \A^{-} \cap \A_r $, for $ r > 0 $.

\bl\label{l8}
The sets $ \A^{-} $, $ \A^{+} $ and $ \A_{r}^{+} $ are non-empty, for all $ r> 0 $.
\el 
\begin{proof}
	Let $ u\in S(c) $. From Remark \ref{obs3}, we can choose $ t_1 >0 $ sufficiently large such that $ u_1 = H(u, t_1) \in S(c) $ and $ V(u_1) < 0 $ and, in a similar way, we can find $ t_2 < 0 $ such that $ u_2 = H(u, t_2) \in S(c) $ and $ V(u_2) \geq 0 $. Moreover, from Remark \ref{obs3}, there exists $ t_3 < 0 $ sufficiently large satisfying $ A(u_3) \leq r $ and $ V(u_3) \geq 0 $, for $ u_3 = H(u, t_3) $. 
\end{proof}

\bl\label{l12}
Let $ \mu \in \R $ and $ (\un) \subset S(c) $ verifying $ \displaystyle \limsup_{n} A(\un) \leq 1 - c $ and $ I(\un) \leq d $, for some $ d\in \R $ and for all $ n \in \N $. Then, there exists a sequence $ (\yn) \subset \mathbb{Z}^2 $ and $u \in X \setminus \{0\}$ such that $ \until=u_n(\cdot-y_n)\rightharpoonup u  $ in $ X $. 
\el 
\begin{proof}
Since $ (A(\un)) $ is bounded, from equation (\ref{v2}), $ (V_2(\un)) $ is also bounded. Moreover, from (\ref{eq2}), (\ref{eq3}) and Corollary \ref{c1}, we see that $ \left(\intR F(\un) dx \right) $ is bounded. Thus, once $ I(\un) \leq d $, for all $ n \in \N $, we have
$$
V_1(\un) = \dfrac{4}{\mu}I(\un) - \dfrac{2}{\mu}A(\un) + V_2(\un) + \dfrac{4}{\mu} \intR F(\un) dx \leq C_1 , \forall \ n \in \N .
$$
From Lemma \ref{l6} and Corollary \ref{c2}, there exists a sequence $ (\yn) \subset \mathbb{Z}^2 $ such that, up to a subsequence, $ \until \CF u $ in $ \Hum $ with $ u \neq 0 $. We can assume, without loss of generality that $ \until(x) \RA u(x) $ a.e. in $ \Rdois $. Thus, once $ V_1(\until) = V_1(\un) $, from Lemma \ref{l2}, $ (||\until||_\ast) $ is bounded. Therefore, up to a subsequence, we conclude that $ \until \CF u $ in $ X $ with $ u \neq 0 $.
\end{proof}

The next lemma is one of the key results to obtain our main theorems.

\bl\label{l13}
Let $ \mu \in \R $ and $ (\un)\subset S(c) $ satisfying $ \displaystyle \limsup_{n} A(\un) \leq 1-c $ and $ I(\un) \leq d $, for some $ d\in \R $ and for all $ n \in \N $. Then, up to a subsequence, $ (\un) $ is bounded in $ X $.
\el
\begin{proof}
First of all, since $ (\un) $ is bounded in $ \Hum $ and from Lemma \ref{l12}, passing to a subsequence if necessary, $ \until=u_n(\cdot-y_n) \CF u $ in $ X $, $ u\neq 0 $ and $ \until(x) \RA u(x) $ pointwise a.e. in $ \Rdois $. Moreover, $ \until \RA u $ in $ \Ldois $. 

\noindent \textbf{Claim:} $(\yn)$ is a bounded sequence. 

Indeed, since $ u\neq 0 $ in $ \Ldois $, there are $ R_1 > 0 $, $ n_1 \in \N $ and $ C_1 > 0 $ such that $ ||\un||_{2, B_{R_1}}^{2}\geq C_1 > 0 $, for all $ n \geq n_1 $. If $|y_n|\leq 2R_1$ for all $n \in \mathbb{N}$ the claim is proved. From this, assume that there is $n \in \mathbb{N}$ such that $ |\yn| > 2R_1 $.  
Recalling that 
$$ 
1+|x+y| \geq 1+\frac{|y|}{2}\geq \sqrt{1+|y|} \quad \forall x \in B_{R_1} \quad \mbox{and} \quad y \in B_{2R_1}^{c},
$$ 
we have that
\begin{align*}
||\until||_{\ast}^{2} =\intR \ln(1+|x|)\un^2 (x-\yn) dx&  = \intR \ln(1+|x+\yn|) \un^2(x) dx \\
& \geq C_2 ||\un||_{2, B_{R_1}}^{2}\ln(1+|\yn|) = C_3 \ln(1+{|\yn|}). 
\end{align*}
As $(||\until||_{\ast})$ is bounded, the above inequality implies that $(y_n)$ is bounded, showing the claim. 

The boundedness of $(y_n)$ yields that $(u_n)$ is bounded in $X$, because $(u_n)$ is bounded in $H^{1}(\mathbb{R}^2)$ and, for all $ n \in \N $, 
$$
\|u_n\|^2_{\ast}=\int_{\mathbb{R}^2}\ln(1+|x|)|u_n|^2\,dx=\int_{\mathbb{R}^2}\ln(1+|x-y_n|)|\tilde{u}_n|^2\,dx \leq \|\tilde{u}_n\|^2_{\ast}+\ln(1+|y_n|)\|\tilde{u}_n\|^2_2 .
$$
The above inequality together with the boundedness of $(\tilde{u}_n)$ in $X$ implies that $ (\un) $ is bounded in $ X $.
\end{proof}

The last result of this section is the other key to obtain our main theorems.

\bl\label{l14}
Let $ \mu \in \R $ and $ (\un) $be a $ (PS) $ sequence for $ I $ restricted to $ S(c) $ bounded in $ X $ and satisfying $\displaystyle \limsup_{n} A(\un) \leq 1-c $. Then, up to a subsequence, $ \un \RA u $ in $ X $ with $ u \neq 0 $. Particularly, $ u $ is a critical point for $ I $ restricted to $ S(c) $.
\el
\begin{proof}
First of all, once $ (\un) $ is bounded in $ X $, passing to a subsequence if necessary, we obtain that $ \un \CF u $ in $ X $ and $ \un \RA u $ in $ \Ls $, for all $ s \geq 2 $. In particular, $ u \in S(c) $. Thus, from Corollary \ref{c1}, we get that $ \left(\intR F(\un) dx \right) $ and $ \left(\intR f(\un) \un dx \right) $ are bounded. Moreover, from (\ref{v2}) and (\ref{eq6}), we obtain that $ (A(\un)), (V_1(\un)) $ and $ (V_2(\un)) $ are also bounded. 

\noindent \textbf{Claim:} There exists a value $ \lambda \in \R $ such that $ (\un) $ is a $ (PS) $ sequence for the functional $ \I(u) \coloneqq I(u) + \frac{\lambda}{2}||u||_{2}^{2} $.

Indeed, once $ (\un) $ is bounded in $ X $, from \cite[Lemma 3]{[lions]}, adapted from the unit sphere to $ S(c) $, we know that $ ||I\big\vert_{S(c)}'(\un)||_{X'}= o(1) $ is equivalent to $ ||I'(\un) - \frac{1}{c}I'(\un)(\un)\un ||_{X'}=o(1) $. 

Set $$ \lambda_n = -\dfrac{1}{c}I'(\un)(\un) = -\dfrac{1}{c}\left[A(\un) + \dfrac{\mu}{4}V(\un)-\intR f(\un) \un dx \right], \forall \ n \in \N . $$
Then, $ (\lambda_n) \subset \R $ is bounded and, up to a subsequence, $ \lambda_n \RA \lambda  $ in $ \R $. Hence, observing that $ \I'(u)=I'(u)+ \lambda u $ and that, for $ v \in X \setminus \{0\} $,
$$
|\I'(\un)(v)| \leq |I'(\un)(v) + \lambda_n \langle \un, v \rangle | + |\lambda - \lambda_n | ||v|| c^{\frac{1}{2}} , \forall \ n \in \N ,
$$
we conclude that $ (\un) $ is a $ (PS) $ sequence for $ \I $, proving the claim.

Now, since $ \un \CF u $ in $ X $ and $ \un \RA u $ in $ \Ls $, for all $ s \geq 2 $, we have

\noindent \textbf{(i)} $ 0 \leq |\I'(u_ n)(\un - u)| \leq ||\I'(\un)||_{X'}||\un - u||_X \RA 0 $.

\noindent \textbf{(ii)} $ \left| \intR f(\un) \un dx \right| \leq \varepsilon ||u_n||_{2\tau}^{\tau}||\un - u||_2 + C_1 ||\un||_{2(q-1)r_2}^{q-1}||\un - u||_{2r_2} \RA 0 $.

\noindent \textbf{(iii)} $ |V_2 '(\un)(\un - u) | \leq C_2 ||\un||_{\frac{8}{3}}^{3}||\un - u ||_{\frac{8}{3}}\RA 0 $.

\noindent \textbf{(iv)} $ V_1 '(\un)(\un -u) = B_1(\un^2 , u(\un - u)) = B_1(\un^2 , (\un - u)^2) + B_1(\un^2 , u(\un - u)) $ .

From Lemma \ref{l3}, $B_1(\un^2 , u(\un - u)) \RA 0  $. Thus, since $ B_1(\un^2 , (\un - u)^2) \geq 0 $ and from (i)-(iv), we obtain 
\begin{align*}
o(1) = \I'(\un)(\un -u) & = o(1) +A(\un)-A(u) + \dfrac{\mu}{4}V'(\un)(\un -u) - \intR f(\un) (\un - u ) dx \\
& \geq o(1) + A(\un) - A(u) \geq o(1), \forall n \in \N .
\end{align*}
Consequently, $ A(\un) \RA A(u) $ and, in particular, $ \un \RA u $ in $ \Hum $. Then, going back to the above inequality, we conclude that $ B_1(\un^2 , (\un - u)^2) \RA 0 $ and, from Lemma \ref{l2}, $ ||\un - u||_\ast \RA 0 $. Therefore, $ \un \RA u $ in $ X $.

Finally, for $ v \in X $, we have
$$
|I\big\vert_{S(c)}'(u)(v)| = \lim | I\big\vert_{S(c)}'(\un) v | \leq ||v||_X \lim ||I\big\vert_{S(c)}'(\un)||_{X'} = 0 ,
$$
which implies that $  u  $ is a critical point to $ I $ restricted to $ S(c) $.
\end{proof}

\section{Proof of Theorems \ref{t1} and \ref{t2}}

Inspired by \cite{[olimpio]}, we will construct a suitable mountain pass level. From Lemma \ref{l7}, for each $ u\in S(c) $ there are a value $ t_u \leq 0 $ such that $ \varphi_u(t) > 0 $ for all $ t \leq t_u $ and a value $ t_{u, A} < 0$ verifying $ A(H(u, t)) < \frac{\rho_c}{2} $ for all $ t \leq t_{u, A} $. Thus, we can define the following real values
$$
- \infty < \overline{t} \coloneqq \sup \{ t_u < 0 \ ; \ u\in \A^{-}  \} \leq 0 \mbox{ \ \ and \ \ }  - \infty < \overline{\overline{t}} \coloneqq \sup \left\{ t_{u, A} < 0 \ ; \ u\in \A^{-} \right\} \leq 0 .
$$
Consider $ t_0 = \min\{\overline{t}, \overline{\overline{t}}\} $. Thus, we can choose $ u_0 \in \A^{-} $ such that there exists $ t_1 < 0 $, with $ t_0 -1 < t_1 < t_0  $, satisfying $ \varphi_{u_0}(t_1)> 0 $ and $ A(H(u_0, t_1)) < \frac{\rho_c}{2} $. Set $ u_1 = H(u_0, t_1) $. Moreover, from Lemma \ref{l7}, there exists $ t_2 > 0 $ such that $ u_2 = H(u_0, t_2) $ satisfies $ A(u_2) > 2\rho_c $ and $ I(u_2) < 0 $. Therefore, we have the mountain pass level $ \m_\nu $ as defined in (\ref{eq10}). 

The above construction makes possible to find a bound to $ \mu $ that does not depend on the fixed function $ u_0 $. As one can observe in the proof of Lemma \ref{l15}, if we simply consider a function $ u_0 $ fixed and a value $ t_1 < 0 $ given by Lemma \ref{l7}, the constant that appears will depend on $ u_0 $. 

\bl\label{l10}
Assume $ (f_1)-(f_3) $ and $ \mu > 0 $ sufficiently small. Then, for $u_1 \in X$ with $A(u_1) \leq \rho_c$, there exists a value $ \rho_c = \rho(c) > 0 $ such that 
$$
0 < I(u_1) < \inf\limits_{u \in \B} I(u) , 
$$
where 
$$
\B = \{ u \in S(c) \ ; \ A(u) = 2 \rho_c \}. 
$$
\el
\begin{proof}
Set $ \rho_c < \frac{1-c}{2} $ and $ v \in S(c)  $ with $ A(v) = 2 \rho_c  $. Then, from equations (\ref{eq2}), (\ref{eq3}), Lemma \ref{l5} and $ v \in S(c) $, we have
$$
\intR F(v) dx \leq C_2 A(v)^{\frac{\tau - 1}{2}} + C_{3}A(v)^{\frac{qr_2 - 2 }{2r_2}} .
$$ 
Thus, since $ F(u_1), V_1(v), V_2(u_1) \geq 0 $ and from (\ref{v2}), we have
\begin{align*}
I(v)- I(u_1) & = \dfrac{1}{2}(A(v)-A(u_1))+ \dfrac{\mu}{4}(V(v) - V(u_1))+ \intR (F(u_1)-F(v)) dx \\
& \geq A(v)-A(u_1)-\dfrac{\mu}{4}V_1(u_1)-\dfrac{\mu}{4}V_2(v)- \intR F(v) dx \\
& \geq \dfrac{1}{2}\rho_c -\dfrac{\mu}{2}\|u_1\|_*^{2}c-\dfrac{\mu}{4}K_1 c^{\frac{3}{2}}\sqrt{2}\rho_{c}^{\frac{1}{2}} - C_4 \rho_{c}^{\frac{\tau - 1}{2}} - C_5 \rho_{c}^{\frac{qr_2 - 2 }{2r_2}} .
\end{align*}
Moreover, once $ c^{\frac{3}{2}}<c $ and $ \rho_{c}^{\frac{1}{2}} < 1 $,  for $ \|u_1\|_*^{2} =C_1 $  we get that
$$
I(v)- I(u_1) \geq \left(\dfrac{1}{4}\rho_c  -\dfrac{\mu}{2}C_1 c  - \frac{\mu}{4}K_1 \sqrt{2} c \right) + \left(\dfrac{1}{4}\rho_c - C_4 \rho_{c}^{\frac{\tau - 1}{2}} - C_5 \rho_{c}^{\frac{qr_2 - 2 }{2r_2}} \right) .
$$
Therefore, fixing $ \rho_c $ even smaller of such way that 
$$
\dfrac{1}{4}\rho_c - C_4 \rho_{c}^{\frac{\tau - 1}{2}} - C_5 \rho_{c}^{\frac{qr_2 - 2 }{2r_2}} \geq \frac{1}{8}\rho_c
$$
and $\mu_c=\mu_c(\rho_c)>0$ such that
$$
\dfrac{1}{4}\rho_c -\dfrac{\mu}{2}C_1 c - \frac{\mu}{4}K_1 \sqrt{2} c \geq \frac{1}{8}\rho_c, \quad \forall \mu \in (0,\mu_c),
$$
we get
$$
I(v)- I(u_1) \geq \frac{1}{4}\rho_c>0,
$$
and so, 
$$ 
\inf\limits_{v\in \B} I(v) > I(u_1).
$$
Moreover, we also have that
$$
I(u_1) \geq \dfrac{1}{2}A(u_1) - \dfrac{\mu}{4}K_1 c^{\frac{3}{2}} A(u_1)^{\frac{1}{2}} - C_1 A(u_1)^{\frac{\tau - 1}{2}} - C_2 A(u_1)^{\frac{qr_2 - 2}{2r_2}} > 0 , 
$$
for $ \rho_c > 0 $ sufficiently small. This proves the desired result.  
\end{proof}

\bl\label{l15}
We have $ \max\limits_{t \in [0, 1]}I(\gamma(t)) > \max\{I(u_1), I(u_2)\} $, for all $ \gamma \in \Gamma $.
\el
\begin{proof}
Let $ \gamma \in \Gamma $. Then, $ ||\gamma(t)||_{2}^{2}=c $ for all $ t \in [0, 1] $. Moreover, $ A(\gamma(0)) < \frac{\rho_c}{2} $ and $ A(\gamma(1)) > 2 \rho_c $. Thus, $ ||\gamma(0)||^2 < c + 2\rho_c < ||\gamma(1)||^2 $ and, from the Intermediate Value Theorem, there exists $ \mathfrak{t}\in [0, 1] $ such that $ ||\gamma(\mathfrak{t})||^2= c + 2 \rho_c $. 

On the other side, $||\gamma(\mathfrak{t})||^2=A(\gamma(\mathfrak{t}))+c$. Consequently, $ A(\gamma(\mathfrak{t})) = 2\rho_c $. Therefore, from Lemma \ref{l10}, $ I(u_1) < I(\gamma(\mathfrak{t})) $ and the result follows.
\end{proof}

As a consequence of Lemma \ref{l15}, $ \m_\nu > 0 $. Next we seek for an useful upper bound for $ \m_\nu $.

\bl\label{l16}
There are $\mu_0, \nu_0>0$ such that $\m_\nu \leq \dfrac{(1-c)(\theta - 6)}{8\theta} $, for $ \nu > \nu_0$ and $ \mu \in (0,\mu_0)$. 
\el
\begin{proof}
Consider the path $ \gamma_0(t)=H(u_0 , (1-t)t_1 + t t_2) \in \Gamma $. Then, 
\begin{align*}
\max\limits_{t \in [0, 1]}I(\gamma_0(t)) & \leq \max\limits_{r\geq 0}\left\{ \dfrac{r}{2}A(u_0) - \dfrac{\mu}{4}c^2(t_0 - 1) - \nu r^{\frac{q-2}{2}}||u_0||_{q}^{q} \right\} \\
& = \max\limits_{r\geq 0}\left\{ \dfrac{r}{2}A(u_0) - \nu r^{\frac{q-2}{2}} ||u_0||_{q}^{q} \right\} + \dfrac{\mu}{4}c^2 K_4 ,
\end{align*}
where $ K_4 = -(t_0 - 1) > 0 $. Hence, we obtain a constant $ K_5 = K_5 (q, u_0) > 0 $ such that
$$
\m_\nu \leq K_5 \left(\dfrac{1}{\nu}\right)^{\frac{2}{q-4}}+K_4 \dfrac{\mu}{4}c^2 ,
$$
and the result follows for 
$$
\mu < \dfrac{(1-c)(\theta - 6)}{4 c^2 \theta K_4}=\mu_0
$$
and 
$$
\nu \geq \left(\frac{16\theta K_5}{(1-c)(\theta-6)} \right)^{\frac{q-4}{2}}=\nu_0.
$$
\end{proof}

As an immediate consequence of the last lemma is the corollary below
\bc \label{obs6}
The Lemma \ref{l16} is true letting $ \mu $ to be any positive real number and controlling the mass $ c $. In this case, it is enough to consider $\nu \geq \nu_0$ and 
\begin{equation} \label{c0}
	c < \min\left\{1, \dfrac{(6-\theta)+ \sqrt{(\theta - 6)^2 + 16 K_4 \mu \theta (6-\theta)}}{8K_4 \mu \theta} \right\}=c_0.
\end{equation}
\ec

In the sequel, let $ (\un)\subset S(c) $ be the sequence constructed in \cite{[jeanjean]} satisfying
\begin{equation}\label{eq11}
I(\un) \RA \m_\nu , \ ||I\big\vert_{S(c)}'(\un)||_{X'}\RA 0 \mbox{ \ \ and \ \ } Q(\un)\RA 0 \mbox{, as } n \RA +\infty , 
\end{equation}
where $ Q: \Hum \RA \R $ is given by
\begin{equation}\label{Q}
Q(u) = A(u) - \dfrac{\mu}{4}||u||_{2}^{4}+2 \intR F(u) dx - \intR f(u) u dx .
\end{equation}

\bl\label{l17}
Let $ (\un) $ be the sequence given in (\ref{eq11}). Then, decreasing if necessary $\mu_0$ given in Lemma \ref{l16}, we have  
$$
\limsup_{n} \intR f(\un)\un dx \leq \dfrac{4\theta}{\theta - 6} \m_\nu, \quad \forall \mu \in (0,\mu_0).
$$
A similar results holds by fixing $\mu$ and decreasing if necessary the number $c_0$ given in (\ref{c0}). 
\el
\begin{proof}
From $ (f_3) $, $ V_1(u) \geq 0 $ and (\ref{v2}), we have 
\begin{equation} \label{COA1}
	\m_\nu + o(1) \geq I(\un) - \dfrac{1}{4}Q(\un) \geq \dfrac{1}{4}A(\un) - \dfrac{\mu}{4}K_1c^{\frac{3}{2}} A(\un)^{\frac{1}{2}}+ \frac{\mu}{16}c^2+ \dfrac{(\theta - 6)}{4 \theta}\intR f(\un) \un dx .
\end{equation} 
Since $ c\in (0, 1) $, it follows that 
$$
\m_\nu + o(1) \geq I(\un) - \dfrac{1}{4}Q(\un) \geq \dfrac{1}{4}A(\un) - \dfrac{\mu}{4}K_1c A(\un)^{\frac{1}{2}}+ \frac{\mu}{16}c^2+ \dfrac{(\theta - 6)}{4 \theta}\intR f(\un) \un dx .
$$
Using Remark \ref{obs5}, with $ \Aa = \frac{1}{4} $, $ \Bb = \frac{\mu}{4}K_1 c  $ and $ \Dd = \frac{\mu}{16}c^2 $, decreasing if necessary $\mu_0$, we obtain
$$
\m_\nu + o(1) \geq \dfrac{(\theta - 6)}{4 \theta}\intR f(\un) \un dx , \forall \ n \in \N .
$$
Now, if we fix $\mu>0$, we apply Remark \ref{obs5}, with $ \Aa = \frac{1}{4} $, $ \Bb = \frac{\mu}{4}K_1 c^{3/2}  $ and $ \Dd = \frac{\mu}{16}c^2 $ in (\ref{COA1}) to get the inequality above,  decreasing if necessary $c_0$ given in (\ref{c0}). 
\end{proof}

\bc\label{c3}
Let $ (\un) $ be the sequence given in (\ref{eq11}). Then, for $\nu \geq \nu_0$ and $\mu >0$, 
$$
\limsup_{n} A(\un) \leq \dfrac{(1-c)}{2} + \dfrac{\mu}{4}c^2 .
$$
\ec
\begin{proof}
From (\ref{eq11}), 
$$
A(\un) = Q(\un) + \dfrac{\mu}{4}c^2 - 2 \intR F(\un) dx + \intR f(\un) \un dx \leq o(1) + \dfrac{\mu}{4}c^2 + \intR f(\un) \un dx .
$$
Now, the result follows employing Lemma \ref{l17}. 
\end{proof}

As an immediate consequence of the last corollary we have 
\bc\label{c4}
Let $ (\un) $ be the sequence given in (\ref{eq11}). Then, decreasing if necessary $\mu_0>0$ given in Lemma \ref{l16} we get $ \displaystyle \limsup_{n} A(\un) \leq 1-c $. Moreover, a similar estimate holds fixing $\mu>0$ and considering  
$$
0<c < \dfrac{-1 + \sqrt{1+ 2 \mu}}{ \mu} .
$$
\ec

\bl\label{c5}
Let $ (\un) $ be the sequence given in (\ref{eq11}). Then, $ (\un) $ is bounded in $ \Hum $.
\el
\begin{proof}
	From Corollary \ref{c4}, there exists $ n_0 \in \N $ such that $ A(\un) \leq 2-c $, for all $ n \geq n_0 $. Thus, $ A(\un) \leq C_2 $, for all $ n \in \N $, where $ C_2 = \max\{ A(u_1),A(u_2),...A(u_{n_0}), 2-c \} $. Therefore, $ ||\un||\leq (C_2 + c )^{\frac{1}{2}} $, for all $ n \in \N $.
\end{proof}

We will only write the proof of Theorem \ref{t1}, because the same argument works to prove Theorem \ref{t2}.  

\noindent {\bf Proof of Theorem \ref{t1}:} Let $ (\un) $ be the sequence given in (\ref{eq11}). Then, from Corollary \ref{c4}, $ \displaystyle \limsup_{n} A(\un) \leq 1-c $. Thus, from Lemma \ref{l14}, without loss of generality, we can assume that $ \un \RA u $  in $ X $, $ u \neq 0 $ and $ u $ is a critical point for $ I $ restricted to $ S(c) $. Moreover, $ I(u)=\m_\nu $.

\section{Proof of Theorems \ref{t3} and \ref{t4}}

This section is devoted to the proof of the multiplicity results. We start recalling the reader the definition of the Krasnoselski's genus. Consider 
$$ 
\Lambda = \{ \K \subset X \ ; \ \K \mbox{ is symmetric and closed} \} \ (\mbox{with respect to topology in \ }  X )
$$ 
and 
$$ 
\Sigma_{\K}= \{ k \in \N \ ; \ \mbox{ there exists } \phi \in C^0 ( \K , \mathbb{R}^k \setminus \{0\} ) \mbox{ such that } \phi(-u) = \phi (u)\}.
$$ 
Then, we define the Krasnoselski's genus of $ \K $, denoted by $ \gamma(\K) $, as follows
$$
\gamma(\K) = \left\{ \begin{array}{ll}
\inf \Sigma_\K  \mbox{ , if \ } \Sigma_\K \neq \emptyset  \\
+ \infty \mbox{ \ \ \ \  , if \ } \Sigma_\K = \emptyset
\end{array}
\right. 
\mbox{ \ \ and \ \ } \gamma(\emptyset) =0 .
$$
Basic properties of genus can be found in \cite[Chapter II.5]{[genus]}. In the sequence, we consider $ Z $ as a real Banach space with norm $ ||\cdot||_Z $,  $ \Hh $ as a Hilbert space with inner product $ \langle \cdot , \cdot \rangle_\Hh $  and $ \M = \{ u \in Z \ ; \ \langle u , u \rangle_\Hh = m \} $, for $ m > 0 $. We assume that $ Z $ is continuously embedded in $ \Hh $. Moreover, define 
$$
\Upsilon_k = \{ \K \subset \M \ ; \ \K \mbox{ \ is simmetric, closed and \ } \gamma(\K) \geq k \} \mbox{ , for \ } n \in \N .
$$
Then, we are ready to enunciate a crucial result, adapted from \cite[Theorem 2.1]{Lu}.

\bt\label{t5}
(\cite[Theorem 2.1]{Ji}) Let $ I : Z \RA \R $ an even functional of $ C^1 $ class. Suppose that $ I\big\vert_\M $ is bounded from below and satisfies the $ (PS)_d $ condition for all $ d< 0 $, and $ \Upsilon_k \neq \emptyset $ for each $ k = 1, ..., n $. Then, the minimax values $ -\infty < d_1 \leq d_2 \leq \cdots \leq d_n  $ can be defined as
$$
d_k = \inf\limits_{\K \in \Upsilon_k}\sup\limits_{u \in \K} I(u) \mbox{ , for \ } 1\leq k \leq n .
$$
Moreover, the following statement are valid:

\noindent \textbf{(i)} $ d_k $ is a critical value for $ I\big\vert_{\M} $, provided $ d_k < 0 $ ; 

\noindent \textbf{(ii)} If $ d_k = d_{k+1} = \cdots = d_{k+l -1} = d < 0 $, for some $ k, l \geq 1 $, then $ \gamma(\K_d) \geq l $, where $ \K_d $ is the set of critical points of $ I\big\vert_\M $ in the level $ d\in \R $. Hence, $ I\big\vert_\M $ has at least $ n $ critical points.
\et  

Our approach is based in \cite{Ji} (see also \cite{Azozero}), which we refer the reader for more details in the arguments. Precisely, in what follows we will verify the necessary conditions to apply Theorem \ref{t5} in order to obtain at least $ n $ solutions for \eqref{P}.

First of all, observe that
$$
I(u) \geq \dfrac{1}{2}A(u) - \dfrac{\mu}{4}K_1 c^{\frac{3}{2}}A(u)^{\frac{1}{2}} - C_1 A(u)^{\frac{\tau - 1}{2}} - C_2 A(u)^{\frac{qr_2 -2}{2r_2}} \geq h(A(u)^{\frac{1}{2}}) , 
$$
for $ C_1, C_2 > 0 $ and $ h: \R \RA \R $ given by
$$
h(t) = \dfrac{1}{2}t^2 - \dfrac{\mu}{4}K_1 c^{\frac{3}{2}} t - C_1 t^{\tau - 1} - C_2 t^{\frac{qr_2 -2}{r_2}} .
$$
Since $ \tau - 1 > 2 $ and $ \frac{qr_2 -2}{r_2} > 2 $, there exists a value $ a > 0 $ such that, if $ \mu c^{\frac{3}{2}} < a $, then there are $ R_0 , R_1 > 0 $ satisfying
$$
\left\{ \begin{array}{lll}
h(t) \leq 0 \mbox{ \ , for \ } t\in [0, R_0],  \\
h(t) \geq 0 \mbox{ \ , for \ } t\in [R_0 , R_1], \\
h(t) < 0 \mbox{ \ , for \ } t\in (R_1 , + \infty).
\end{array}
\right. 
$$ 
Moreover, we can define the values
\begin{equation}\label{eqe1}
\mu_1 = \dfrac{a}{c^{\frac{3}{2}}} \mbox{ \ \ \  and \ \ \ } c_1 = \left(\dfrac{a}{\mu}\right)^{\frac{3}{2}}. 
\end{equation}

Now, for $ R_0 , R_1 $ given above, define $ \T : \R^+ \RA [0, 1] $ as a non-decreasing function such that $ \T \in C^{\infty} $ and
$$
\T(t) = \left\{ \begin{array}{ll}
1 \mbox{ , for \ } t \in [0, R_0], \\
0 \mbox{ , for \ } t \in [R_1 , + \infty).
\end{array}
\right.
$$
Thus, we consider the truncated functional $ I_\T : X \RA \R $ given by
$$
I_\T (u) = \dfrac{1}{2}A(u) + \dfrac{\mu}{4} V(u) - \T (A(u)^{\frac{1}{2}}) \intR F(u) dx .
$$
Similarly as above, we have 
\begin{equation}\label{eqe2}
I(u) \geq \overline{h}(A(u)^{\frac{1}{2}}) , 
\end{equation}
where $ \overline{h}: \R \RA \R $ is defined as
$$
\overline{h}(t) = \dfrac{1}{2}t^2 - \dfrac{\mu}{4}K_1 c^{\frac{3}{2}} t - \T(t) [ C_1 t^{\tau - 1} + C_2 t^{\frac{qr_2 -2}{r_2}} ] .
$$
Without loss of generality, we assume that 
$$
\dfrac{1}{2}t^2  - C_1 t^{\tau - 1} - C_2 t^{\frac{qr_2 -2}{r_2}} \geq 0 \ , \ \forall t \in [0, R_0 ], \mbox{ \ \ and \ \ } R_0 < \sqrt{1-c} .
$$

\bl\label{le1}
\begin{itemize}
\item[(1)] $ I_\T \in C^1 (X, \R) $.

\item[(2)] If $ I_\T (u) \leq 0 $, then $ A(u)^{\frac{1}{2}} < R_0 $ and $ I(v) = I_\T (v) $ for all $ v $ in a small neighbourhood of $ u $ in $ X $. 

\item[(3)] If $ \mu \in (0, \mu_1) $, then $ I_\T $ restricted to $ S(c) $ verifies the $ (PS)_d $ condition in every level $ d < 0 $. 
\end{itemize}
\el
\begin{proof}
Items (1) and (2) can be proved by standard arguments, so we omit it here. Let us prove item (3). Let $ (\un) \subset S(c) $ be a $ (PS) $ sequence for $ I_\T $ restricted to $ S(c) $ in the level $ d < 0 $. Then, up to a subsequence, $ I(\un) \leq 0 $, for all $ n \in \N $. Thus, from item (2), we conclude that $ \displaystyle \limsup_{n} A(\un) \leq 1-c $. Consequently, the result follows from Lemmas \ref{l13} and \ref{l14}. 
\end{proof}

In what follows we will need the level sets
$$
I_{\T}^{d} = \{ u \in S(c) \ ; \ I_\T (u) \leq d \} . 
$$

\bl\label{le2}
For each $ n \in \N $ and $ \mu \in \R $, there are $ \epsilon_n = \epsilon(n) > 0 $ and $ \nu_n = \nu(n) > 0 $ such that $ \gamma(I_{\T}^{-\epsilon}) \geq n $, for all $ \epsilon \in (0, \epsilon_n) $ and $ \nu \geq \nu_n $.
\el
\begin{proof}
For each $ n \in \N $, as in \cite{Ji}, consider the $ n$-dimensional space $ E_n \subset X $ with the orthogonal base $ \B = \{ u_1 , ... , u_n \} $, that is, 
$$
\intR \nabla u_j \nabla u_k dx =\intR u_j u_k dx = \intR \ln (1+|x|) u_j u_k dx=0 ,
$$
if $ j\neq k $, $ A(u_j) = \rho^2 < 1-c $, $ ||u_j||_{2}^{2}=c $ and $ ||u_j|| = \sqrt{\rho^2 + c} $, for each $ 1\leq j \leq n $. Define 
$$
Z_n = \{ t_1 u_1 + \cdots t_n u_n \ ; \ t_{1}^{2} + \cdots t_{n}^{2} = 1 \} \mbox{ \ \ and \ \ } \Ss_{\rho^2 + c} = \{ (y_1, ... , y_n ) \in \mathbb{R}^n \ ; \  y_{1}^{2}+ \cdots y_{n}^{2}= \rho^2 + c \} .
$$
Considering the map $ \Phi : Z_n \RA \Ss_{\rho^2 + c} $ given by $ \Phi(u) = (t_1 \sqrt{\rho^2 + c}, ... , t_n \sqrt{\rho^2 + c} ) $, for $ u = \sum\limits_{j=1}^{n}t_j u_j $, one can easily see that $ Z_n $ and $ \Ss_{\rho^2 + c} $ are homeomorphic (with respect to $ X $). Thus, by genus properties, we get that $ \gamma(Z_n) = n $. 

Now, since $ \dim E_n < + \infty $, all the norms are equivalent. Then, the value
$$
a_n = \inf \left\{ ||u||_{q}^{q} \ ; \ u \in S\left(\frac{c}{\rho^2}\right) \cap E_n \mbox{ \ and \ } A(u) = 1 \right\}
$$
is well-defined and positive. Moreover, there exists a constant $ P_n = P(n) > 0 $ such that $ ||u||_{\ast}^{2} \leq P_n A(u) $, for all $ u \in E_n $. 

Observe that, once $ \B $ is orthogonal, $ A(v)=\rho^2 $ for all $ v \in Z_n $ and, considering $ 0 < \rho < R_0 $, we have $ I(v)=I_\T(v) $. Hence, from condition $ (f_4) $,  we have
\begin{align*}
I_\T (v) = I(v) & \leq \dfrac{1}{2}\rho^2 A\left(\dfrac{v}{\rho}\right) + \dfrac{\mu}{4}\rho^4 \left\Vert \dfrac{v}{\rho} \right\Vert_{2}^{2}\left\Vert \dfrac{v}{\rho} \right\Vert_{\ast}^{2} - \nu \rho^q \left\Vert \dfrac{v}{\rho} \right\Vert_{q}^{q} \\
& \leq  \dfrac{1}{2}\rho^2 + \dfrac{\mu}{4}P_n c \rho^2 - \nu a_n \rho^q , \ \forall \ v \in Z_n . 
\end{align*}
Therefore, choosing $ \rho \in (0, R_0) $ there is $ \epsilon_n > 0 $ and a value $ \nu_n > 0 $ sufficiently large such that $ I_\T(v) \leq - \epsilon_n $, for all $ v\in Z_n $ and $\nu \geq \nu_n$. Consequently, $ Z_n \subset I_{\T}^{- \epsilon} $ and, from a genus property, $ \gamma(I_{\T}^{-\epsilon}) \geq n $. 
\end{proof}

\bl\label{le3}
Let $ \Gamma_k =\{ D \subset S(c) \ ; \ D \mbox{ is symmetric, closed and \ } \gamma(D) \geq k \} $, 
$$
d_k = \inf_{D \in \Gamma_k}\sup\limits_{u\in D} I_\T (u) 
$$
and $ \K_d = \{ u\in S(c) \ ; \ I_\T '(u) = 0 \mbox{ \ and \ } I_\T (u) = d\}   $. Assume that $ \mu \in (0, \mu_1) $, where $ \mu_1 $ is given by \eqref{eqe1}. If $ d_k < 0 $ then $ d_k $ is a critical value of $ I_\T $ restricted to $ S(c) $. Moreover, if $ d_k = \cdots d_{k+r} \coloneqq d < 0 $, for some $ k, r \geq 1 $, then $ \K_d \neq \emptyset $ and $ \gamma(\K_d) \geq r+1 $. Particularly, $ I_\T $ restricted to $ S(c) $ has at least $ k $ non-trivial critical points.
\el
\begin{proof}
From Lemma \ref{le2}, for each $ k \in \N $, there exists $ \epsilon_k > 0 $ such that $ \gamma(I_{\T}^{-\epsilon}) \geq k $, for all $\epsilon \in (0, \epsilon_k)$. Since $ I_\T $ is continuous and even, $ I_{\T}^{-\epsilon} \in \Gamma_k $. Consequently, $ d_k \leq - \epsilon_k < 0 $, for all $ k \in \N $. On the other hand, from \eqref{eqe2}, $ I_\T $ is bounded from bellow over $ S(c) $, which implies $ d_k > - \infty $, for all $ k \in \N $. Hence, from Lemma \ref{le1}-(3) and Theorem \ref{t5}, $ d_k $ is a critical value of $ I_\T $ restricted to $ S(c) $. 

In the sequence, suppose that $ d_k = \cdots d_{k+r} \coloneqq d < 0 $, for some $ k, r \geq 1 $. From the first part, $ \K_d \neq \emptyset $. Observe, now, that a sequence $ (\un) \subset \K_d $ is a $ (PS)_d $ sequence. Hence, from Lemma \ref{le1}-(1), $ R_0 < \sqrt{1-c} $ and Lemma \ref{l14}, $ \K_d $ is compact. Finally, from a deformation lemma, genus properties and Theorem \ref{t5}, the result follows.
\end{proof}

\begin{proof}[Proof of Theorems \ref{t3} and \ref{t4}]
One should observe that, from Lemma \ref{le1}-(2), critical points to $ I_\T $ restricted to $ S(c) $ are precisely critical points of $ I $ restricted to $ S(c) $. Therefore, the proof follows from Lemma \ref{le3}.
\end{proof}
%\section{References}

\vspace{\baselineskip} 

\noindent \textbf{Acknowledgements:} C.O. Alves was supported by  CNPq/Brazil 304804/2017- ; E.S. Böer was supported by Coordination of Superior Level Staff Improvement-(CAPES)-Finance Code 001 and  S\~ao Paulo Research Foundation-(FAPESP), grant $\sharp $ 2019/22531-4, and O.H. Miyagaki was supported by  National Council for Scientific and Technological Development-(CNPq),  grant $\sharp $ 307061/2018-3 and FAPESP  grant $\sharp $ 2019/24901-3. 

\begin{center}{\bf Declarations} \end{center}
{\bf Conflict of Interest.} On behalf of all authors, the corresponding author states that there is no conflict of
interest. \\
{\bf Data Availability Statement.} This article has no additional data.

\end{document}